\documentclass[11pt,a4paper]{article}
\usepackage{amssymb,latexsym,amsmath,amsthm,hyperref} 

\def\a{\alpha} 
\def\b{\beta} 
\def\l{\lambda} 
 
\def\u{\mu} 
\def\ka{\kappa} 
\def\th{\theta} 
 
\def\r{\rho} 
\def\e{\eta} 
\def\x{\xi} 
\def\z{\zeta}

\def\sg{\sigma} 
\def\Sg{\Sigma}

\def\R{\mathbb R} 
\def\C{\mathbb C} 
\def\H{\mathbb H} 
\def\S{\mathbb S} 
\def\P{\mathbb P} 
\def\12{\frac{1}{2}} 
\def\dd{{\rm d}} 

\def\zb{\bar z} 
\def\pa{\partial} 
\def\ti{\tilde}
\def\wti{\widetilde} 
\def\<{\langle} 
\def\>{\rangle} 
\providecommand{\abs}[1]{\lvert#1\rvert} 

\newtheorem{theorem}{Theorem}[section]
\newtheorem{lemma}[theorem]{Lemma} 
\newtheorem{corollary}[theorem]{Corollary} 
\newtheorem{proposition}[theorem]{Proposition}

\theoremstyle{remark} 
\newtheorem{remark}[theorem]{Remark} 
\newtheorem{example}[theorem]{Example} 

\theoremstyle{definition} 
\newtheorem{definition}[theorem]{Definition} 

\title{Adjoint Transform of Willmore Surfaces in $n$-sphere}
\author{Xiang Ma}

\begin{document}

\maketitle

\begin{abstract}
After the surface theory of M\"obius geometry, this study 
concerns a pair of conformally immersed surfaces in $n$-sphere. 
Two new invariants $\th$ and $\r$ associated with them 
are introduced as well as the notion of touch and co-touch. 
This approach is helpful in research about transforms of 
certain surface classes. As an application, we define 
adjoint transform for any given Willmore surface in $n$-sphere. 
It always exists locally (yet not unique in general) 
and generalizes known duality theorems of Willmore surfaces. 
This theory on surface pairs reaches its high point by a 
characterization of adjoint Willmore surfaces in terms of 
harmonic maps.
\end{abstract}

\section{Introduction}

One fascinating aspect of surface theory in differential geometry is
the construction of various transforms preserving certain surface
classes \cite{quater, Jeromin, Rogers+Schief}. Many of them
are classical and discovered more than one hundred years ago. 
Today geometers are still interested in such transforms, because they
indicate a hidden symmetry for the surface class in concern, which shows
the deeper connection with integrable systems \cite{Burstall, Cieslinski, 
Jeromin, Rogers+Schief}. Conversely, if a surface class is associated 
with some integrable equation or enjoys similar features, 
one would expect to find such transforms.

In this work we are interested in construction of new transforms for 
Willmore surfaces in $\S^n$. It has been shown that they allow a 
\emph{spectral transform} similar to isothermic surfaces \cite{BPP}. 
Both surface classes are M\"obius
invariant and connected by the classical Blaschke's Problem 
(\cite{Blaschke},\cite[Ch.3]{Jeromin}). This problem requires to find
two surfaces enveloping the same 2-sphere congruence and forming 
conformal correspondence. The non-trivial solutions in $\S^3$ consist 
of Darboux pair of isothermic surfaces and dual Willmore surfaces. 
Since there already exists a rich transform theory of isothermic surfaces
\cite{Burstall, Jeromin}, we are encouraged to find a parallel theory
of Willmore surfaces.

This aim was partially achieved in \cite{quater}. For a Willmore surface 
in $\H P^1\cong \S^4$ they introduced three kinds of transforms: 
the forward and backward 1-step B\"acklund transform, which resemble 
the Christoffel transform of an isothermic surface (also known as the 
dual isothermic surface); the forward and backward 2-step B\"acklund 
transform; the Darboux transform, which is described by a Riccati equation
like the Darboux transform of isothermic surfaces \cite{Jeromin+Pedit}.

Although the quaternionic setup provides new insight into the spinor 
representation of surfaces in 3- and 4-space, it does not apply
to higher dimensional spaces. For this reason it is not realistic to
generalize these transforms to $\S^n$ based on the original algebraic 
description. Therefore, we follow another line: \emph{For a pair of 
surfaces being certain transform to each other, characterize them 
by geometric conditions}. One of such results is the
well-known geometric characterization of B\"acklund transforms for
pseudo-spherical surfaces \cite[Ch.6]{Terng}. Another example is
Blaschke's Problem and its solutions mentioned above. In view of 
such results, it is reasonable to establish a general theory on 
surface pairs at the beginning. That should be done in a M\"obius 
invariant way in arbitrary dimensional space. We accomplished 
this task in the light-cone model, and derived two new invariants 
$\th,\r$ associated with such a pair of oriented surfaces \cite{Ma}. 

Another difficulty of establishing the transform theory of Willmore 
surfaces lies in the following fact: there exists no \emph{dual surface}
for a generic Willmore surface in $\S^n$. Notice that the dual 
isothermic surface is the basic transform for a given isothermic 
surface, on which the construction of all other transforms rely 
\cite{Burstall, Jeromin+Pedit}. So the failure of Bryant's duality 
theorem (\cite{Bryant}, see also \cite{Blaschke}) in higher 
codimension case is really a disappointment. 

Thus it might be a surprise to the reader that there does exist a 
generalization of dual Willmore surface to $\S^n$, which is 
introduced in this paper and called the \emph{adjoint Willmore 
surface(s)}.
We used the plural at here, because in general they are not unique.
Yet such \emph{adjoint transforms} always exist locally, including 
the dual S-Willmore surface as a special case. We show the adjoint 
transform always produces a Willmore surface $\hat{f}$ from a given 
one $f:M\to \S^n$, and $f$ is also an adjoint transform of $\hat{f}$
(Theorem~\ref{thm-duality}).
This generalizes the known duality theorems of Bryant \cite{Bryant}, 
Ejiri \cite{Ejiri} and that about the forward/backward 2-step 
B\"acklund transforms \cite{quater}. Indeed, the definition of the 
adjoint Willmore surface is inspired by the geometric characterization
of the 2-step B\"acklund transform in terms of \emph{co-touching
condition}, the case in which our new invariant $\th=0$. 

In another viewpoint, a surface pair is just a (Riemann) surface mapped 
into the moduli space of point pairs in $\S^n$. This map is conformal
exactly when $\th=0$. Furthermore, if we ask what is a 
conformal harmonic map into this semi-Riemannian symmetric space,
we find it must be given by a pair of adjoint Willmore surfaces
(Theorem~\ref{thm-harmonic}). This interesting result 
enjoys a similar flavor as the well-known characterization
of Willmore surfaces by the harmonicity of its conformal Gauss map
(see Theorem~\ref{thm-willmore}).

We make two remarks relating our new results to the Darboux transform
of isothermic surfaces. First, the condition of \emph{touching} 
($\r=0$), the counterpart of co-touching, was originally introduced 
in \cite{Bohle}, subsequently used to define the generalized Darboux 
transform for a generic surface in $\S^4$, the moduli space of which
gives the \emph{spectral curve}. Next, it is worth 
mentioning that the symmetric space of point pairs has been 
considered for $\S^3$ \cite{Burstall-flat}, $\S^4$ \cite{Jeromin-flat}
and $\S^n$ \cite{Burstall}\cite[Ch.8]{Jeromin}. 
In these works, Darboux pair of isothermic surfaces is characterized 
in terms of \emph{curved flats} in this symmetric space, which is a
special kind of integrable system. Again we find some parallel between
the theory of Willmore surfaces and isothermic surfaces.

In the following, we will briefly review the surface theory in 
M\"obius differential geometry in Section~\ref{sec-surface} and some
basic facts about Willmore surfaces. We go on to develop the theory
on surface pairs in Section~\ref{sec-pair} with a detailed discussion
on the meaning of touch and co-touch (the interpretation by quaternions
is left to the Appendix). Finally, the definition of adjoint transform
of a Willmore surface as well as its properties are presented in 
Section~\ref{sec-adjoint}.


\section{Surface theory in M\"obius geometry} 
\label{sec-surface} 

\subsection{The surface theory by moving frames} 

In this paper we will follow \cite{BPP} in their treatment of 
surface theory in M\"obius geometry. As usual, let $\mathcal{L}$ 
denote the 
light cone in the $n+2$ dimensional Minkowski space $\R^{n+1,1}$ 
with quadratic form $\< y,y\>=-y_0^2+\sum_{i=1}^{n+1}y_i^2$. 
Then the unit sphere $\S^n(n\ge 3)$ in Euclidean space may be identified 
with our projectivized light cone: 
\[ 
\S^n \cong \P(\mathcal{L}):x\leftrightarrow [1:x]. 
\] 
The projective action of the Lorentz group on $\P(\mathcal{L})$ 
yields a representation of the M\"obius group. 
In this model, points are described by light-like vectors 
(null lines), and hyperspheres correspond to space-like vectors. 
Generally, a $k$-sphere $S\subset\S^n$ is represented by space-like 
$(n-k)$-dim subspace $U$ (or the orthogonal complement $U^{\bot}$
equivalently). 

For surface $f:M\to \S^n \cong \P(\mathcal{L})$, 
a \emph{(local) lift} of $f$ is just a map $F$ from $M$ 
into the light cone such that the null line spanned by 
$F(p)$ is $f(p)$. Two different local lifts differ by a scaling, 
so the metric induced from them are conformal to each other.
When the underlying $M$ is a Riemann surface, $f$ is a conformal 
map iff $\< F_z,F_z\>=0$ for any $F$ and any coordinate $z$ on $M$; 
it is immersion iff $\< F_z,F_{\zb}\> > 0$. 

In our study we often associate a 2-sphere (congruence) to $f$ 
and be interested in the contact relationship. Let space-like 
$(n-2)$-dim subspace $U$ stand for such a $2$-sphere, $F$ a lift of $f$. 
This sphere passes through $f(p)$ iff $\< F(p),U\>=0$. 
Suppose this is satisfied, then the sphere is tangent to $f$ 
at $p$ iff $\< \dd F(p),U\>=0.$ Identify the 2-sphere 
with $U^{\bot}$, then it is tangent to $f$ at $p$ iff $F(p)$ 
(the map itself) and $\dd F(p)$ (all tangent vectors) are
contained in $U^{\bot}$. 

Given conformal immersion $f:M\to \S^n \cong \P(\mathcal{L})$ 
of Riemann surface $M$ with local lift $F$, there is a M\"obius 
invariant decomposition 
$M\times\R^{n+1,1} = V \oplus V^{\bot}$, 
where 
\[ 
V=\mbox{Span}\{F,\dd F,F_{z\zb}\} 
\] 
is a rank-$4$ subbundle defined via local lift $F$ and complex 
coordinate $z$ (one readily checks that $V$ is independent to 
such choices, thus well-defined). $V$ is a Lorentzian subbundle, 
and $V^{\bot}$ is a space-like subbundle, which might be identified 
with the normal bundle of $f$ in $\S^n$. The connection $D$ on $V^{\bot}$ 
defined by orthogonal projection of the derivative in $\R^{n+1,1}$ 
is the usual \emph{normal connection} in metric geometry, which is 
already known to be M\"obius invariant. On the other hand, 
$V$ determines a M\"obius invariant 2-sphere 
$\mathbb{P}(V\cap\mathcal{L})$ 
at every point of this immersed surface. we call it the 
\emph{mean curvature sphere} or \emph{central sphere congruence}.
The complexification of $V$ and $V^{\bot}$ are denoted respectively 
as $V_{\C},V_{\C}^{\bot}.$ 

\begin{remark} 
\label{rem-mean} 
The name \emph{mean curvature sphere} comes from the remarkable property 
that it is tangent to the surface and has the same mean curvature vector 
as the surface at the tangent point, where the ambient space is endowed 
with a metric of Euclidean space (or any space form).
\end{remark} 

Fix a local coordinate $z$. Among various choice of local lifts there is 
a canonical one into the forward light cone, which is denoted by $Y$.
$Y$ is M\"obius invariant and determined by 
$\abs{\dd Y}^2=\abs{\dd z}^2. $ We choose a M{\"o}bius invariant frame 
of $V_{\C}$ as $ \{Y,Y_z,Y_{\zb},N\}$.
The real $N\in \Gamma(V)$ is chosen so that these frame vectors are 
orthogonal to each other except 
$ 
\< Y_z,Y_{\zb} \> = \12,~\< Y,N \> = -1. 
$ 
Such a light-like vector $N$ is also unique. 

Since $Y_{zz}$ is orthogonal to $Y$, $Y_z$ and $Y_{\zb}$, there must be 
a complex function $s$ and a section $\ka\in\Gamma(V_{\C}^{\bot})$ 
so that the following Hill's equation holds: 
\begin{equation} 
\label{eq-schwarz} 
Y_{zz} + \frac{s}{2} Y = \ka. 
\end{equation} 
This defines two basic invariants $\ka$ and $s$ depending on 
coordinate $z$. $\ka$ may be identified with the normal valued 
\emph{Hopf differential} up to a suitable scaling, meanwhile 
$s$ is interpreted as the \emph{Schwarzian} of immersion $f$.
In $\S^3$, $\ka$ and $s$ form a complete system of invariants.
For more explanation, see \cite{BPP}.

Let $\psi\in\Gamma(V^{\bot})$ denote an arbitrary section of the 
normal bundle. Now it is easy to derive the structure equations: 
\begin{equation}\label{eq-structure} 
\left\{ 
\begin{aligned} 
   Y_{zz} &= -\frac{s}{2} Y + \ka, \\[-0.1cm]
   Y_{z\zb} &= -\<\ka,\bar{\ka}\> Y + \12 N, \\ 
   N_{z} &= -2 \<\ka,\bar{\ka}\> Y_z - s Y_{\zb}+ 2D_{\zb}\ka, \\[.1cm]
   \psi_z &= D_z\psi+2\<\psi,D_{\zb}\ka\> Y-2\<\psi,\ka\> Y_{\zb}. 
\end{aligned} 
\right. 
\end{equation} 
The computation is straightforward, hence omitted at here. The conformal 
Gauss, Codazzi and Ricci equations as integrable conditions 
are given as below: 
\begin{subequations} 
\begin{gather} 
  \12s_{\zb} = 3\< D_z\bar{\ka},\ka\> + \<\bar{\ka},D_z\ka\>, 
   \label{gauss}\\ 
  \mbox{Im}\Big( D_{\zb}D_{\zb}\ka + \frac{\bar{s}}{2}\ka\Big) = 0, 
   \label{codazzi}\\ 
  R_{\zb z}^D\psi := D_{\zb}D_z\psi - D_z D_{\zb}\psi 
   = 2\<\psi,\ka\>\bar{\ka}- 2\<\psi,\bar{\ka}\>\ka. 
   \label{ricci} 
\end{gather} 
\end{subequations}

\subsection{Willmore functional and Willmore surfaces} 
\label{subsec-Willmore} 

\begin{definition} 
For a conformally immersed surface $f:M\to\S^n$ with decomposition 
$M\times\R^{n+1,1} = V \oplus V^{\bot}$ as before, we define 
\[ 
G:=Y\wedge Y_u\wedge Y_v\wedge N
=-2i\cdot Y\wedge Y_z\wedge Y_{\zb}\wedge N,~\quad z=u+iv. 
\] 
It is a map from $M$ to the Grassmannian $G_{3,1}(\R^{n+1,1})$,
called the \emph{conformal Gauss map} of $f$. This Grassmannian 
consists of all 4-dimensional Minkowski subspaces. A basic fact
about this map is \cite{Ejiri}
\end{definition} 
\begin{proposition} 
\label{prop-conformal}
For conformal immersion $f:M\to\S^n$, $G$ induces a positive 
definite metric (by the usual inner product between multivectors)
\[ 
g=\frac{1}{4}\<\dd G,\dd G\>=\<\ka,\bar\ka\>\abs{\dd z}^2 
\] 
on $M$ except at umbilic points, which is called the 
\emph{M\"obius metric}.
Especially this is a conformal metric, 
thus justifies the name of conformal Gauss map. 
\end{proposition} 
\begin{definition} 
The \emph{Willmore functional} of $f$ is defined at here as 
the area of $M$ with respect to the M\"obius metric: 
\[ 
W(f):=\frac{i}{2}\int_M \abs{\ka}^2 \dd z\wedge\dd\zb. 
\footnote{
In case of a surface in $\R^3$, the Willmore functional is usually
defined as $\widetilde{W}(f)=\int_M (H^2-K)\dd M. $ 
It differs from our definition by $\widetilde{W}(f)=4W(f)$.
} 
\] 
\end{definition}

\begin{definition} 
Let $M$ be a topological surface. Any immersion $f:M\to\S^n$ 
automatically induces a conformal structure over $M$, hence 
defines the Willmore functional $W(f)$. If $f$ is a critical point 
of $W$ with respect to any variations of the map and the induced 
conformal structures, it is called a \emph{Willmore surface}. 
\end{definition} 

For any conformal map $G:M \to G_{3,1}(\R^{n+1,1})$, the energy is
$E(G):=\int_M \<\dd G\wedge*\dd G\>$.
The Willmore functional of a surface $f$ is related to 
the energy of its conformal Gauss map via
$W(f)=-\frac{1}{8}E(G).$
Moreover, Willmore surfaces are characterized by 
the harmonicity of its conformal Gauss map. 
\begin{theorem} 
[\cite{Bryant,BPP,Ejiri}] 
\label{thm-willmore} 
For a conformally immersed surface $f$ in $\S^n$, 
the following three conditions are equivalent: 
\begin{enumerate} 
\item[(i)] 
  $f$ is Willmore. 
\item[(ii)] 
  The Hopf differential and Schwarzian of $f$ satisfy 
\begin{equation} 
\label{willmore1} 
D_{\zb}D_{\zb}\ka + \12\bar{s}\ka = 0. 
\quad\text{(Willmore condition)} 
\end{equation} 
  This is a condition stronger than the conformal Codazzi 
equation~\eqref{codazzi}. 
\item[(iii)] 
  The conformal Gauss map $G$ is a harmonic map into 
the Grassmannian $G_{3,1}(\R^{n+1,1})$. 
\end{enumerate} 
\end{theorem} 
\begin{corollary}
\label{cor-spectral}
The integrability conditions for a Willmore surface is
\begin{equation*} 
\left\{
\begin{gathered}
\12 s_{\zb} = 3\< D_z\bar{\ka},\ka\> + \<\bar{\ka},D_z\ka\>, \\ 
D_{\zb}D_{\zb}\ka + \frac{\bar{s}}{2}\ka = 0, \\ 
R_{\zb z}^D\psi= 2\<\psi,\ka\>\bar{\ka}- 2\<\psi,\bar{\ka}\>\ka. 
\end{gathered}
\right.
\end{equation*} 
This system admits the symmetry
\[ 
\ka_{\l}=\l\ka,\quad s_{\l}=s, 
\]
for unitary $\l\in S^1$, which describes the associated family of 
Willmore surfaces.
\end{corollary}
\begin{remark}
\label{rem-spectral}
The characterization of Willmore surfaces in terms of \eqref{willmore1} 
was given in \cite{BPP} without proof. Note the equivalence 
between conditions (ii) and (iii) of Theorem~\ref{thm-willmore} 
is well-known to experts in this field, from which (i) follows easily. 
Analogous results in Lie sphere geometry and projective geometry
were discussed in \cite{Burstall+Jeromin}.
\end{remark}

\section{Pair of conformally immersed surfaces} 
\label{sec-pair} 

After the general surface theory, we turn to transforms for 
certain surfaces. Usually they are obtained from a given surface 
by some integrable equations. Alternatively, many times such 
transforms might be characterized by some geometric conditions 
on the surface pair involved. The second approach motivates us to
build a general theory of surface pairs, which seems to be a 
natural development based on the previous section.

\subsection{Basic invariants of surface pair} 
\label{subsec-pair} 

Let us start with a Riemann surface $M$ and two conformal immersions 
$f,\hat{f}:M\to \S^n$ which are assumed to be always distinct. 
Given coordinate $z$, set $Y$ to be the canonical lift of $f$, 
with Schwarzian $s$, Hopf differential $\ka$, 
and frame $\{Y,Y_z,Y_{\zb},N\}$. $\widehat Y$ is
a fixed local lift of $\hat{f}$ so that $\< Y,\widehat{Y}\>=-1. $ 
We may express explicitly that   
$\widehat{Y}=\l Y + \bar{\u}Y_z + {\u}Y_{\zb} + N + \x,$ 
where $\l$ and $\u$ are real-valued and complex-valued functions 
respectively, and the real $\x\in\Gamma(V^{\bot}).$ 
Since $\widehat{Y}$ is isotropic, there must be 
$\l=\12(\abs{\u}^2+\<\x,\x\>).$ So we have
\begin{equation} 
\label{yhat1} 
\widehat{Y}=\12\big(\abs{\u}^2+\<\x,\x\>\big) Y 
+ \bar{\u}Y_z + {\u}Y_{\zb} + N + \x, 
\end{equation} 
Take derivative on both sides. By \eqref{eq-structure} we may find the
fundamental equation for such a surface pair after a straightforward 
computation: 
\begin{equation} 
\label{yhat-z1} 
\widehat{Y}_z =\frac{\u}{2}\widehat{Y} 
+ \th\left(Y_{\zb} + \frac{\bar\u}{2}Y\right) 
+ \r\left(Y_z + \frac{\u}{2}Y\right) + \<\x,\z\> Y + \z , 
\end{equation} 
where 
\begin{subequations} 
\label{eq-theta-rho} 
\begin{align} 
\th &= \u_z - \12{\u}^2 - s - 2\<\x,\ka\>, \label{theta1} \\ 
  \r &= \bar{\u}_z - 2\<\ka,\bar\ka\> 
       + \12\<\x,\x\>, \label{rho1} \\ 
  \z &= D_z\x - \frac{\u}{2}\x + 
       2\left(D_{\zb}\ka + \frac{\bar\u}{2}\ka\right) 
       ~\in~\Gamma(V_{\C}^{\bot}). 
       \label{zeta1} 
\end{align} 
\end{subequations} 
It is easy to check that $\th$ and $\r$ corresponds to a $(2,0)$ form 
and a $(1,1)$ form separately. They may be defined alternatively
by the inner product between bivectors: 
\begin{subequations} 
\label{thetarho-geom}
\begin{align}
\th &= 2~\< Y\wedge Y_z,\widehat{Y}\wedge\widehat{Y}_z\>, 
            \label{theta-geom}\\ 
\r &= 2~\< Y\wedge Y_{\zb},\widehat{Y}\wedge\widehat{Y}_z\>. 
            \label{rho-geom} 
\end{align} 
\end{subequations} 
Although these expressions involve lifts $Y,\widehat{Y}$ and 
coordinate $z$, both $\th$ and $\r$ are independent to such choices, hence
well-defined invariants associated with such a pair of immersed surfaces.
Note if we interchange between $Y$ and $\widehat{Y}$, $\r$ turns to be 
$\bar\r$, and $\th$ keeps invariant. \\

$Y\wedge Y_z$ and $\widehat{Y}\wedge\widehat{Y}_z$ may be 
interpreted as \emph{complex contact elements}. So
$\th$ and $\r$ are second order invariants of $(f,\hat{f})$.
More concretely, a 2-dim contact element (always assumed to be oriented) 
is just a 2-dim oriented subspace in $T_p\S^n$ for some $p\in\S^n$, 
which corresponds to a 3-dim oriented subspace of signature $(2,0)$ 
in $\R^{n+1,1}$. We represent this object by its oriented frame 
$\{X,X_1,X_2\}$ with scalar product matrix $\mbox{diag}(0,1,1)$. 
This determines (up to multiplication by an unitary complex number) the 
\emph{complex contact element} represented by $X\wedge(X_1-{\rm i}X_2)$. 
Its conjugate corresponds to the real contact element with reversed 
orientation. 

Consider two contact elements $\Sg=\{Y,Y_1,Y_2\}$ and 
$\widehat\Sg=\{\widehat Y, \widehat Y_1,\widehat Y_2\}$
at distinct points (so $\< Y,\widehat{Y}\>\ne 0$). 
Similarly define
\begin{subequations} 
\begin{align} 
\th
&= \12\frac{\< Y\wedge(Y_1-{\rm i}Y_2), 
   \widehat{Y}\wedge(\widehat{Y}_1-{\rm i}\widehat{Y}_2)\>} 
   {\< Y\wedge\widehat{Y},Y\wedge\widehat{Y}\>}, 
   \label{theta-contact}\\ 
\r
&= \12\frac{\< Y\wedge(Y_1+{\rm i}Y_2), 
   \widehat{Y}\wedge(\widehat{Y}_1-{\rm i}\widehat{Y}_2)\>} 
   {\< Y\wedge\widehat{Y},Y\wedge\widehat{Y}\>}. 
   \label{rho-contact} 
\end{align} 
\end{subequations} 
Note they are independent to the choice of frames.

When $\< Y,\widehat{Y}\>=0$, we have two contact elements at the same point. 
Intuitively we need only to consider the 2-planes $\mbox{Span}\{Y_1,Y_2\}$ and 
$\mbox{Span}\{\widehat Y_1,\widehat Y_2\}$. The following two quantities 
\begin{subequations} 
\label{theta-rho}
\begin{align} 
\underline\th 
&=\12\< Y_1+{\rm i}Y_2,\widehat Y_1-{\rm i}\widehat Y_2\>,\\ 
\underline\r 
&=\12\< Y_1-{\rm i}Y_2,\widehat Y_1-{\rm i}\widehat Y_2\>. 
\end{align} 
\end{subequations} 
are similarly well-defined, i.e. they are independent to the 
choice of frames of $\Sg,\widehat\Sg$. 
Compared to \eqref{theta-contact} \eqref{rho-contact}, here 
the $\pm$ sign is reversed in two places. 
Why this convention will be clear in next subsection.

\subsection{Touch and co-touch} 
\label{subsec-touch} 

To better understand the geometric meaning of $\th$ and $\r$
(as well as their counterparts $\underline\th,\underline\r$), 
let's consider the special case when either of them vanishes. 

\begin{definition} 
Two contact elements $\Sg$ and $\widehat\Sg$ at one point 
are said to \emph{touch} each other if $\underline\r=0$ and 
\emph{co-touch} each other if $\underline\th=0$. 

Consider two oriented surfaces immersed in $\S^n$ intersecting at $p$. 
We say they \emph{touch (co-touch)} each other if the contact elements 
given by their tangent spaces at $p$ touch (co-touch). 
\end{definition} 

\begin{example} 
\label{exa-tangent} 
For two surfaces tangent to each other at the same point, 
it is easy to see that they either touch each other 
at this point when their orientations are compatible, 
or co-touch when the orientations are opposite. 
\end{example} 

\begin{example} 
\label{exa-complex} 
Given two complex lines in $\C^n, n\ge 2.$ Regard them as 
real 2-planes with the induced orientation (via the complex structure) 
in $\R^{2n}$. Then they touch each other. In Appendix we will 
see that all touching 2-plane pairs are constructed in this way. 
\end{example} 

These examples show that the touching relation (including touch
and co-touch) between two surfaces is a generalization of tangency.
Such notions were first introduced by Pedit and Pinkall 
in the context of quaternions $\H$, then used to define 
\emph{Darboux transforms} for general surfaces in $\S^4\cong \H P^1$,
which generalize the classical Darboux transforms of 
isothermic surfaces \cite[Section~7.1]{Bohle}. 
Simply speaking, for a surface immersed in $\H\cong \R^4$ 
one can define the left and right normal vector $N,R$. 
Two surfaces having a common point $p$ are said to \emph{left
touch} each other if they share the same left normal vector $N$ at $p$.
\footnote{
In other words, left-touching means the tangent planes of these
surfaces at $p$ can be transformed to each other by 
right-multiplication of a unit quaternion. Right touch is 
understood in the similar way.
}
Co-touching was defined similarly \cite{Ma0} when they have 
opposite $N$ or $R$. Detailed discussion is left to Appendix.

The key observation in \cite{Ma0} is: the touching relation 
between 2-planes depends only on the Euclidean geometry 
of $\R^4$, whereby independent to the quaternionic structure. 
Moreover, it might be defined in arbitrary dimensional space. 
This seems out of one's expectation and calls for explanation. 
Note the usual Gauss map identifies a 2-plane in $\R^n$ with 
a null line in $Q\subset \C P^{n-1}$. For two such null lines
$l_1,l_2\in Q$, there are two noteworthy cases, i.e. when $l_1$
is orthogonal to $l_2$, or to $l_2$'s conjugate. Either of these
two cases corresponds to touch or co-touch between the 
2-planes represented by $l_1,l_2$. 
\footnote{
For a pair of intersecting lines in $\R^n$, the intersection angle
is the only invariant in Euclidean geometry.
More generally, for two oriented $m$-dim subspaces 
$\Sg_1,\Sg_2\subset \R^n$, the singular values of the inner product 
matrix between their oriented orthonormal frames are a complete system 
of invariants under the action of $\rm{SO}(n)$.
They are independent to the choice of such frames, hence well-defined.
In M\"obius geometry we find the complete invariants associated with
a pair of oriented contact elements at the same point in this way. 
When $m=2$, suppose the singular values are $\l_1,\l_2$, 
then $\Sg_1$ touch (co-touch) $\Sg_2$ iff $\l_1=\l_2$ ($\l_1=-\l_2$).
}\\

To clarify the geometric meaning of $\th=0$ and $\r=0$ for a pair 
of conformal immersions with lifts $Y, \widehat Y$, observe that given 
coordinate $z=u+{\rm i}v$, contact element $\Sg=\{Y,Y_u,Y_v\}$ 
at $Y(p)$, and single point $\widehat Y(p)$, there is an unique oriented 
2-sphere passing through $Y(p),\widehat Y(p)$ and tangent to $Y$ 
with compatible orientation. It is given by the 4-dim subspace of 
signature $(3,1)$ spanned by $\{Y,Y_u,Y_v,\widehat Y\}$, with the orientation 
fixed by the oriented contact element $\Sg=\{Y,Y_u,Y_v\}$ or the 
complexification $Y\wedge Y_z$. Denote it as $S(p)$. Now we may state 

\begin{proposition} 
\label{prop-touch} 
Given two conformal immersions $f,\hat{f}$, the invariant $\r(p)=0$ iff 
the 2-sphere $S(p)$ touches $\hat{f}$ at $\widehat Y(p)$, and $\th(p)=0$ iff 
$S(p)$ co-touches $\widehat Y$ at $\hat{f}(p)$. 
\end{proposition} 
\begin{proof} 
We may take the normalized lifts $Y,\widehat Y$ as before. By \eqref{yhat1}, 
\[ 
\widehat{Y}=\12\biggl(\abs{\u}^2+\<\x,\x\>\biggr)Y+\bar\u Y_z+\u Y_{\zb}+N+\x 
\] 
is orthogonal to $Y_z + \frac{\mu}{2}Y$. Note that under the reflection 
with respect to $Y-\widehat Y$, $S(p)$ is invariant with reversed orientation, 
and the complex contact element $\Sg=Y\wedge(Y_z+\frac{\mu}{2}Y)$ 
is mapped to $\widehat Y\wedge(Y_z+\frac{\mu}{2}Y)$. Thus the complex contact 
element given by $S(p)=\mbox{Span}\{Y,Y_u,Y_v,\widehat Y\}$ at $\widehat Y(p)$ 
should be $\Sg'=\widehat Y\wedge(Y_{\bar z}+\frac{\bar\mu}{2}Y)$. 
On the other hand, the complex contact element given by immersion 
$\widehat Y$ at $\widehat Y(p)$ is $\widehat\Sg=\widehat Y\wedge \widehat Y_z$. 
Thus at $\widehat Y(p)$ the invariants associated with $\Sg'$ and $\widehat\Sg$ 
are computed by the fundamental equation \eqref{yhat-z1}: 
\[ 
\underline\th=2\< Y_z+\frac{\mu}{2}Y,\widehat Y_z\>=\th,~~ 
\underline\r=2\< Y_{\bar z}+\frac{\bar\mu}{2}Y,\widehat Y_z\>=\r. 
\] 
The conclusion now follows from the definition of touch and co-touch. 
\end{proof}

\section{Adjoint transforms of Willmore surfaces} 
\label{sec-adjoint} 

\subsection{Motivation and definition} 
\label{subsec-definition} 

After Bryant's work \cite{Bryant}, people are interested in the 
generalization of the duality theorem for Willmore surfaces 
in $\S^n$. Ejiri pointed out that the duality theorem 
holds true only for a smaller class of Willmore surfaces, 
the so-called \emph{S-Willmore} surfaces \cite{Ejiri}.
Although that, the hope to generalize the construction of dual 
Willmore surface still exists, according to the observations below:
\begin{enumerate}
\item 
In \cite{Ma}, as an application of our theory on surface pairs, 
Blaschke's Problem and its solutions were generalized to $\S^n$.
Dual S-Willmore surfaces arises as the second class of non-trivial
solutions, for which the invariant $\th=0$ (in the isothermic 
case, $\r=0$). The vanishing of $\th$ has a nice geometric 
interpretation as \emph{co-touching} in general case.
\item
The forward and backward 2-step B\"acklund transforms of a 
Willmore surface in $\S^4$ \cite{quater} are generalization 
of the duality theorem above, which might be called the 
\emph{left and right dual} Willmore surface respectively. 
Like the dual Willmore surface in 3-dim case, the 2-step B\"acklund 
transform also falls on the mean curvature sphere $S$ of the given
Willmore surface. But it only co-touches $S$ and not necessarily 
to be tangent.
\end{enumerate}
Stimulated by these facts, one naturally attempts to characterize
the left and right dual Willmore surface in $\S^4$ by the geometric
properties listed above. It yields a new class of transforms for any
Willmore surface in $\S^n$.
\begin{definition} 
A map $\hat{f}:M\to\S^n$ is called the adjoint transform of Willmore 
surface $f:M\to\S^n$ if it is conformal and co-touches 
the mean curvature sphere of $f$ at corresponding point. Especially, 
$\hat{f}$ must locate on the corresponding mean curvature sphere of $f$. 
Note that $\hat{f}$ is allowed to be a degenerate point. 
\end{definition} 
This definition gives the conditions characterizing an adjoint transform. 
Yet we need a more explicit description. Consider surface pair $f,\hat{f}$ 
with adapted lifts $Y,\widehat{Y}$, satisfying $\< Y,\widehat{Y}\>=-1.$ 
Furthermore suppose $\hat{f}$ is on the mean curvature sphere of $f$. 
Then equations \eqref{yhat1}\eqref{yhat-z1} take the form 
\begin{align} 
\widehat{Y} 
   &=\12\abs{\u}^2 Y + \bar\u Y_z + \u Y_{\zb} + N, 
   \label{yhat}\\ 
\widehat{Y}_z 
   &=\frac{\u}{2}\widehat{Y} 
   + \th\left(Y_{\zb} + \frac{\bar{\u}}{2}Y\right) 
   + \r\left(Y_z + \frac{\u}{2}Y\right) + 2\e. 
   \label{yhat-z} 
\end{align} 
Here $\u\dd z$ is a complex connection 1-form determined by 
$\mu=2\<\widehat{Y},Y_z\>$. It further defines those invariants associated 
with the pair $f,\hat{f}$ as in Subsection~\ref{subsec-pair}: 
\begin{equation}  
\th := \u_z - \12{\u}^2 - s,~~ 
\r := \bar{\u}_z - 2\<\ka,\bar{\ka}\>, ~~
\e := D_{\zb}\ka + \frac{\bar\u}{2}\ka.   
\end{equation} 
There follows 
\begin{alignat}{2} 
\<\widehat{Y}_z,\widehat{Y}_z\> 
&=\<\widehat{Y}_z-\frac{\u}{2}\widehat{Y}, 
\widehat{Y}_z-\frac{\u}{2}\widehat{Y}\> 
&&=4\<\e,\e\>+\th\cdot\r~, \label{yzyz}\\ 
\<\widehat{Y}_z,\widehat{Y}_{\zb}\> 
&=\<\widehat{Y}_z-\frac{\u}{2}\widehat{Y}, 
\widehat{Y}_{\zb}-\frac{\bar{\u}}{2}\widehat{Y}\> 
&&=4\<\e,\bar{\e}\>+\12\abs{\th}^2+\12\abs{\r}^2~. 
\label{yzyzbar} 
\end{alignat} 
That $f$ is Willmore implies 
\begin{equation} 
\label{willmore} 
0=D_{\zb}D_{\zb}\ka+\frac{\bar{s}}{2}\ka 
=D_{\zb}(\e-\frac{\bar{\u}}{2}\ka)+\frac{\bar{s}}{2}\ka 
=D_{\zb}\e-\frac{\bar{\u}}{2}\e-\frac{\bar{\th}}{2}\ka 
\end{equation} 
\begin{definition} 
The map into $\S^n$ represented by \eqref{yhat} is an 
\emph{adjoint transform} of Willmore surface $Y$ 
iff $\u$ satisfies the following conditions: 
\begin{subequations} 
\label{eq-adjoint} 
\begin{alignat}{2} 
&\text{\bf Co-touching:} 
&\quad& 0=\th=\u_z - \12{\u}^2 - s.\label{coto} \\ 
&\text{\bf Conformality:} 
&\quad& 0=\<\e,\e\>=\< D_{\zb}\ka + \frac{\bar\u}{2}\ka, 
D_{\zb}\ka + \frac{\bar\u}{2}\ka\>. \label{cofo} 
\end{alignat} 
\end{subequations} 
\end{definition} 
\begin{example} 
\label{exa-adjoint1} 
A Willmore surface $f$ is a S-Willmore surface if $D_{\zb}\ka$ 
linearly depends on $\ka$. In such a case there exist a function 
$\u$ locally so that $D_{\zb}\ka + \frac{\bar\u}{2}\ka=0$ when
there is no umbilic points. It is easy to check that 
\eqref{eq-adjoint} holds for this $\u$, which gives the dual 
Willmore surface $\hat{f}$ via \eqref{yhat}
\end{example} 
\begin{remark}
It is easy to show that $\u\dd z$ is a connection 1-form of $K^{-1}$, 
where $K$ denotes the canonical bundle of Riemann surface $M$.
Conversely, given any connection 1-form of $K^{-1}$, if it satisfies
\eqref{eq-adjoint} with respect to any local coordinate, then it defines
an adjoint transform globally. 
\end{remark}

\begin{remark}
The reader should be aware of the problem of singularities. 
First, the map underlying $\widehat{Y}$ may not be immersion when $\sg=0$. 
Thus $\wti{Y}$ as well as the underlying map $\ti{f}:M\to\S^n$ 
might has branch points. Next, the connection 1-form 
$\u\dd z$ may have \emph{poles}, which corresponds to the coincidence 
case of $f$ and $\hat{f}$. In this paper we will concentrate on the
local aspect of this construction, and ignore this problem temporarily.
But when deal with closed Willmore surfaces, this is an inevitable 
problem related to both global and local geometry. 

In \cite{Ma0}, this adjoint transform is applied to the study of 
Willmore 2-spheres in $\S^n$, which yields very strong vanishing results. 
Despite this success as well as the removable singularity theorem utilised
there, these singularities still constitute the final obstruction to a
complete classification. If we can have a better understanding of them,  
the known classification results of Willmore 2-spheres \cite{Bryant,
Ejiri, Montiel} might be generalized to $\S^n$ based on \cite{Ma0}.
\end{remark}

\subsection{Existence} 
\label{subsec-existence} 

Our definition of adjoint transforms leads to the natural problem of 
existence and uniqueness of solutions to system \eqref{coto}\eqref{cofo}. 
Note that when $\<\ka,\ka\> \neq 0$, \eqref{cofo} is a 
quadratic equation about $\u$ and much easier to solve. 
In such a situation, at every point we have two roots for 
\[ 
0=\<\e,\e\>=\< D_{\zb}\ka + \frac{\bar\u}{2}\ka, 
D_{\zb}\ka + \frac{\bar\u}{2}\ka\>. 
\] 
Fix either of such a root $\u$ and differentiate this equation. 
By \eqref{willmore}, 
\[ 
0=\<\e,\e\>_{\zb}=2\< D_{\zb}\e,\e\> 
=2\< \frac{\bar\u}{2}\e+\frac{\bar\th}{2}\ka,\e\> 
=\bar\th\<\ka,\e\>. 
\] 
If $\<\ka,\e\>\neq 0$, we have $\th=0$ as desired. Otherwise, 
suppose $\<\ka,\e\>=0$ on an open subset and take derivative, 
one obtains 
\begin{multline*} 
\qquad\qquad 
0=\<\ka,\e\>_{\zb} 
=\< D_{\zb}\ka,\e\>+\<\ka,D_{\zb}\e\>\\ 
=\< \e-\frac{\bar\u}{2}\ka,\e\> 
+\<\ka,\frac{\bar\u}{2}\e+\frac{\bar\th}{2}\ka\> 
=\frac{\bar\th}{2}\<\ka,\ka\>. 
\qquad\qquad 
\end{multline*} 
By assumption, $\<\ka,\ka\> \neq 0$, so $\th=0$. Hence we see that 
the Willmore condition \eqref{willmore} guarantees a solution $\u$ 
of \eqref{eq-adjoint} and the existence of adjoint transforms. 
\\ 
\\ 
How about the case when $\<\ka,\ka\>=0$ on an open subset? 
By Willmore condition \eqref{willmore1} it follows 
$0=\< D_{\zb}\ka,\ka\>=\< D_{\zb}\ka,D_{\zb}\ka\>.$ 
That means \eqref{cofo} holds automatically for any $\u$. 
So we need only to solve the PDE \eqref{coto} 
\[ 
\u_z-\12{\u}^2-s=0 
\] 
independently. It is a Riccati equation about $\u$ with respect to 
the given Schwarzian $s$. 
In S-Willmore case this is solved in Example~\ref{exa-adjoint1}. 
When immersion $f:M\to\S^n$ is Willmore but not S-Willmore, 
the Willmore condition \eqref{willmore1} implies that $\ka$ 
and $D_{\zb}\ka$ span a rank 2 holomorphic subbundle of $V_{\C}^{\bot}$ 
\footnote{ 
This is true at least on the open subset where 
$\ka\wedge D_{\zb}\ka\ne 0$. Then this subbundle extends smoothly
to zeros of $\ka\wedge D_{\zb}\ka$. Note $V_{\C}^{\bot}$ 
is a holomorphic bundle w.r.t. $D_{\zb}$. 
}. 
By \eqref{willmore}, it is easy to show that there is a 1-1 
correspondence between solution $\u$ and holomorphic line subbundle 
spanned by $D_{\zb}\ka+\frac{\bar\u}{2}\ka$.
So there are infinitely many solutions $\u$. 

\begin{remark} 
Generally speaking, \eqref{coto} is a under-determined equation, 
thus admit (infinitely) many solutions. More concretely, there is a 
well known correspondence between the solutions of Riccati equation 
$\u_z-\12{\u}^2-s=0$ and the solutions to linear equation
\begin{equation}
\label{eq-second}
y_{zz}+\frac{s}{2}y=0.
\end{equation}
Suppose of 
\eqref{eq-second} and they are independent, i.e. there exist no
anti-holomorphic function $h$ so that $\ti{y}=h\cdot y$. 
The general solution to \eqref{eq-second} is a combination 
$hy+\ti{h}\ti{y}$, where $y,\ti{y}$ are two independent non-trivial 
solutions, $h,\ti{h}$ are all anti-holomorphic.
In case that only one non-trivial solution $y$ is given, the second 
solution $\ti{y}$ might be found by solving a $\pa$-problem for 
$\l$:~$\l_z=1/y^2$, then $\ti{y}=\l y$. This implies that even for 
a S-Willmore surface with $\<\ka,\ka\>\equiv 0$, locally there are 
still infinitely many adjoint transforms. 
\end{remark}

\subsection{Duality theorem} 
\label{subsec-duality} 

In this subsection, we want to prove that the adjoint transform
preserve the Willmore condition.
Fix the original Willmore surface with lift $Y$. Assume 
there is a $\u$ solving \eqref{coto} and \eqref{cofo}, 
which defines an adjoint transform $\hat{f}$. 
Therefore, \eqref{yhat-z} is simplified to 
\begin{equation} 
\label{yhat-z-new} 
\widehat{Y}_z=\frac{\u}{2}\widehat{Y} 
+ \r\left(Y_z+\frac{\u}{2}Y\right) + 2\e. 
\end{equation} 
Note $\th=0$ also implies 
\begin{equation} 
\label{rho-zbar}
\r_{\zb} 
=\bar\u_{\zb z}-2\<\ka,\bar\ka\>_{\zb} 
=\bar{s}_z + \bar\u\bar\u_z-2\<\ka,\bar\ka\>_{\zb} 
=\bar\u\r +4\<\e,\bar\ka\> 
\end{equation} 
by Gauss equation \eqref{gauss}, and
\begin{equation}
\label{willmore2}
D_{\zb}\e =\frac{\bar\u}{2}\e 
\end{equation}
by \eqref{willmore}. Next consider the canonical lift of 
the adjoint transform, denoted as $\wti{Y}$. 
Let $\<\wti{Y},Y\>=-1/\sg$ be a real function defined 
on $M$. Equivalently speaking, $\wti{Y}$ is obtained from $\widehat{Y}$ 
via $\wti{Y}=\frac{1}{\sg}\widehat{Y}.$ So 
$\12=\<\wti{Y}_z,\wti{Y}_{\zb}\> 
=\frac{1}{{\sg}^2}\<\widehat{Y}_z,\widehat{Y}_{\zb}\>. $
Combined with \eqref{yzyzbar} and $\th=0$, we get 
\begin{equation} 
\label{sigma2} 
{\sg}^2=2\<\widehat{Y}_z,\widehat{Y}_{\zb}\> 
=8\<\e,\bar{\e}\>+\abs{\r}^2. 
\end{equation} 
Thus \eqref{yhat-z-new} may be written as 
\begin{equation} 
\label{ytilde-z} 
\wti{Y}_z=-\frac{\ti\u}{2}\wti{Y} 
+\frac{\r}{\sg}\left(Y_z + \frac{\u}{2}Y\right)+\frac{2}{\sg}\e, 
\end{equation} 
where 
\begin{equation} 
\label{mutilde} 
\ti\u:= \frac{2\sg_z}{\sg} - \u. 
\end{equation} 
To find $\wti{N}$ we should calculate $\wti{Y}_{z\zb}$. It is easy to find 
\begin{subequations} 
\begin{align} 
\e_{\zb} 
&=2\<\e,\bar{\e}\> Y 
-2\<\e,\bar{\ka}\>\left(Y_z+\frac{\u}{2}Y\right) 
+\frac{\bar{\u}}{2}\e,\label{eta-zbar}\\ 
\e_z 
&=-2\<\e,\ka\>\left(Y_{\zb}+\frac{\bar{\u}}{2}Y\right) 
+D_z\e.\label{eta-z} 
\end{align} 
\end{subequations} 
Differentiate both sides of \eqref{ytilde-z}. 
A straitforward computation yields 
\[
\wti{Y}_{z\zb}
= \12\big(\r-\ti\u_{\zb}\big)\wti{Y} 
+\12\Big(-\12\abs{\ti\u}^2\wti{Y} 
-\ti\u\wti{Y}_{\zb}-\bar{\ti\u}\wti{Y}_z+\sg Y\Big). 
\]
Define 
\begin{equation}\label{ntilde} 
\wti{N}:=-\12\abs{\ti\u}^2\wti{Y} 
-\ti\u\wti{Y}_{\zb}-\bar{\ti\u}\wti{Y}_z+\sg Y. 
\end{equation} 
We verify $ \<\wti{N},\wti{Y}_z\>=0,~\<\wti{N},\wti{Y}\>  =-1,~
\<\wti{N},\wti{N}\> = 0. $
So $\{\wti{Y},\wti{Y}_z,\wti{Y}_{\zb},\wti{N}\}$ 
is the canonical frame of $\wti{Y}$ as desired. 
Compare the structure equation (of $\wti{Y}$) 
\[
\wti{Y}_{z\zb} =-\<\ti\ka,\bar{\ti\ka}\>\wti{Y} +\12\wti{N}
\] 
with previous result, we may similarly define 
\begin{equation}\label{rhotilde} 
\ti\r:=\bar{\ti\u}_z 
-2\<\ti\ka,\bar{\ti\ka}\>. 
\end{equation} 
Then there must be 
\begin{equation}\label{rhorho} 
\ti\r=\bar{\r}. 
\end{equation} 

How about the corresponding invariants $\ti\ka$ and $\ti{s}$~? 
According to structure equations of $\wti{Y}$, $\ti{s}$ is determined by 
$\ti\ka=\wti{Y}_{zz}+\frac{\ti{s}}{2}\wti{Y}\in\ti{V}^{\bot}$, 
where $\ti{V} 
:=\mbox{Span}\{\wti{Y},\wti{Y}_z,\wti{Y}_{\zb},\wti{Y}_{z\zb}\} 
=\mbox{Span}\{\wti{Y},\wti{Y}_z,\wti{Y}_{\zb},Y\}$ by \eqref{ntilde}. 
Since $\wti{Y}_{zz}$ and $\wti{Y}$ 
are always orthogonal to $\wti{Y},\wti{Y}_z,\wti{Y}_{\zb}$, 
We find $\ti{s}$ by solving
\[
0=\<\wti{Y}_{zz}+\frac{\ti{s}}{2}\wti{Y},Y\> 
= \<\wti{Y}_z,Y\>_z-\<\wti{Y}_z,Y_z\> -\frac{\ti{s}}{2\sg}
= \frac{1}{2\sg} 
\left(\ti\u_z-\12\ti\u^2-\ti{s}\right). 
\]
Therefore 
\begin{equation} 
\label{stilde} 
\ti{s}=\ti\u_z-\12\ti\u^2. 
\end{equation} 
Denote the normal connection of $\wti{Y}$ as $\wti{D}$. 
We have structure equation 
\begin{equation}
\label{eq-diff-tildekappa} 
2\wti{D}_{\zb}\ti\ka
=\wti{N}_z +2\<\ti\ka,\bar{\ti\ka}\>\wti{Y}_z +\ti{s}\wti{Y}_{\zb}. 
\end{equation}  
Differentiate \eqref{eq-diff-tildekappa} 
and modulo components of $\wti{V}$, which is spanned by $\{\wti{Y},\wti{Y}_z,\wti{Y}_{\zb},Y\}$, one obtains 
\begin{align*} 
2\big(\wti{D}_{\zb}\ti\ka\big)_{\zb} 
+\bar{\ti{s}}\ti\ka 
&\equiv \wti{N}_{z\zb}+\ti{s}\wti{Y}_{\zb\zb} 
+\bar{\ti{s}}\ti\ka \\ 
&\equiv \left(-\12\abs{\ti\u}^2\wti{Y} 
-\ti\u\wti{Y}_{\zb}-\bar{\ti\u}\wti{Y}_z 
+\sg Y\right)_{z\zb}+\ti{s}\wti{Y}_{\zb\zb} 
+\bar{\ti{s}}\ti\ka \\ 
&\equiv \big(\ti{s}-\ti\u_z\big)\wti{Y}_{\zb\zb} 
-\bar{\ti\u}_{\zb}\wti{Y}_{zz} 
-\ti\u\wti{Y}_{\zb\zb z} 
-\bar{\ti\u}\wti{Y}_{zz\zb} 
+(\sg Y)_{z\zb}+\bar{\ti{s}}\ti\ka \\ 
&\equiv \big(\ti{s}-\ti\u_z\big)\bar{\ti\ka} 
+\big(\bar{\ti{s}}-\bar{\ti\u}_z\big){\ti\ka} 
-\ti\u\bar{\ti\ka}_z 
-\bar{\ti\u}\ti\ka_{\zb} 
+(\sg Y)_{z\zb} \\ 
&\equiv -\12\ti\u^2\bar{\ti\ka} 
-\12\bar{\ti\u}^2{\ti\ka} 
-\ti\u\wti{D}_z\bar{\ti\ka} 
-\bar{\ti\u}\wti{D}_{\zb}\ti\ka 
+(\sg Y)_{z\zb} \\ 
&\equiv -\12\ti\u^2\bar{\ti\ka} 
-\12\bar{\ti\u}^2{\ti\ka} 
-\frac{\ti\u}{2}\wti{N}_{\zb} 
-\frac{\bar{\ti\u}}{2}\wti{N}_z 
+(\sg Y)_{z\zb} \\ 
&\equiv -\12\ti\u^2\bar{\ti\ka} 
-\12\bar{\ti\u}^2{\ti\ka} 
-\frac{\ti\u}{2}\left[-\ti\u\wti{Y}_{\zb\zb} 
+(\sg Y)_{\zb}\right] \\ 
&\quad\qquad-\frac{\bar{\ti\u}}{2} 
\left[-\bar{\ti\u}\wti{Y}_{zz} 
+(\sg Y)_z\right]+(\sg Y)_{z\zb} \\ 
&\equiv -\frac{\ti\u}{2}(\sg Y)_{\zb} 
-\frac{\bar{\ti\u}}{2}(\sg Y)_z 
+(\sg Y)_{z\zb} \\ 
&\equiv \left(\sg_z-\frac{\ti\u}{2}\sg\right)Y_{\zb} 
+\left(\sg_{\zb}-\frac{\bar{\ti\u}}{2}\sg\right)Y_z 
+\sg\cdot\12N \\ 
&\equiv \frac{\sg}{2} 
\big({\u}Y_{\zb}+\bar{\u}Y_z+N\big) \\ 
&\equiv \frac{\sg}{2}\left(\sg\wti{Y} 
-\12\abs{\u}^2 Y\right) \\ 
&\equiv 0. 
\end{align*}
Thus we have proved that the Willmore condition \eqref{willmore} is also 
satisfied for $\wti{Y}$. Furthermore, equation \eqref{ntilde} shows 
that $Y$ may be viewed as an adjoint transform of $\wti{Y}$, 
because $\ti\u$ satisfies \eqref{stilde}, which amounts to say 
that the similarly defined quantity $\ti{\th}$ also vanishes, meanwhile 
we already know the \emph{conformality} between $\wti{Y}$ and $Y$.  
This remarkable \emph{duality} is just what one expected, since 
such a relationship between a Willmore surface in $\S^4$ and its 
\emph{forward/backward two-step B{\"a}cklund transforms}(\cite{quater})
is already known. Sum together, we get 

\begin{theorem}
\label{thm-duality}
An adjoint transform $\wti{Y}$ of a Willmore surface $Y$ is also Willmore, 
which is called an \emph{adjoint Willmore surface of $Y$} or \emph{a Willmore 
surface adjoint to $Y$}. Vice versa, $Y$ is also an adjoint transform of 
$\wti{Y}$. The relationship between their corresponding invariants are 
given by  \[ 
-\frac{1}{\sg}=\<\wti{Y},Y\>,~~ 
\frac{2\sg_z}{\sg}=\ti\u+\u,~~ 
\ti\r=\bar{\r}. \] 
\end{theorem} 

\subsection{Characterization by conformal harmonic maps} 
\label{subsec-harmonic} 

In last section, we have developed a theory of pairs of conformally
immersed surfaces $f,\hat{f}$ from a Riemann surface $M$ into $\S^n$. 
Another way to look at them is considering the 2-plane spanned by 
their lifts $Y,\widehat{Y}$. This defines a map
\begin{align*}
H:M &\to G_{1,1}(\R^{n+1,1}),\\
p &\mapsto Y(p)\wedge\widehat{Y}(p).
\end{align*}
Similar to the description of the conformal Gauss map, here the 
Grassmannian $G_{1,1}(\R^{n+1,1})$ consists of all 2-dim Minkowski 
subspaces, and we regard it as a submanifold embedded in 
$\wedge^2\R^{n+1,1}$. The bivector is uniquely determined 
if we put the restriction $\< Y,\widehat{Y}\>=-1$ (hence $\< H,H\>=-1$).
Conversely, such a map corresponds to a pair of surfaces in $\S^n$. 

Associated with $Y,\widehat{Y}$ are invariants $\th,\r$ defined via 
\eqref{eq-theta-rho}. They appear also as invariants of $H$. 
It turns out
\footnote
{In the computation, without loss of generality we may 
assume $Y$ is the canonical lift of $f$, and using the formulae in 
Subsection~\ref{subsec-pair}.}
\begin{gather}
\< H_z,H_z\>=\th,\quad
\< H_z,H_{\zb}\>=\12(\r+\bar\r). \notag\\
\Longrightarrow\quad
\< \dd H,\dd H\>=\th\dd z^2+\12(\r+\bar\r)
(\dd z\dd\zb+\dd\zb\dd z)+\bar\th\dd\zb^2.  \label{eq-H}
\end{gather}
So the co-touching condition is equivalent to the conformality of $H$. 

On the other hand, \eqref{eq-H} gives the energy of $H$:
\[
E(H):=\int_M \< \dd H\wedge *\dd H\>
=-{\rm i}\int_M (\r+\bar\r) \cdot\dd z\wedge\dd\zb
\]
Now comes another natural question: What is the condition that $H$ 
being conformal harmonic? Note that $H$ is similar to the conformal 
Gauss map in that each of them is into some Grassmannian associated 
with $\R^{n+1,1}$. The latter being harmonic iff the original surface 
is Willmore (Theorem~\ref{thm-willmore}). By analogy one would expect 
some similar result for $H$. Of course we should assume that the 
underlying maps $f,\hat{f}$ are also conformal. Surprisingly, these 
simple conditions give a nice characterization of adjoint Willmore surfaces. 
\begin{theorem}
\label{thm-harmonic}
Let $M$ be a Riemann surface. 
Assume $Y,\widehat{Y}$ are local lifts of immersions 
$f,\hat{f}:M\to\S^n$ satisfying $\< Y,\widehat{Y}\>=-1$, 
which induce map 
$H=Y\wedge\widehat{Y}:M\to G_{1,1}(\R^{n+1,1})$. 
Then the three conditions below are equivalent:
\begin{enumerate}
\item[(i)] 
$f,\hat{f}$ are two Willmore surfaces adjoint to each other.
\item[(ii)]
$f,\hat{f}$ and $H$ are conformal maps, and $f,\hat{f}$ 
locate on the mean curvature sphere of each other.
\item[(iii)]
$f,\hat{f}$ are conformal to each other and
$H=Y\wedge\widehat{Y}$ is conformal harmonic. 
\end{enumerate}
\end{theorem}
\begin{proof}
Choose $Y,\widehat{Y}$ as in Subsection~\ref{subsec-pair}, with
\begin{align*} 
\widehat{Y}
&=\12\biggl(\abs{\u}^2+\<\x,\x\>\biggr) Y 
+ \bar{\u}Y_z + {\u}Y_{\zb} + N + \x, \\ 
\widehat{Y}_z 
&=\frac{\u}{2}\widehat{Y} 
+ \th\left(Y_{\zb} + \frac{\bar\u}{2}Y\right) 
+ \r\left(Y_z + \frac{\u}{2}Y\right) + \<\x,\z\> Y + \z, 
\end{align*}  
where $\th,\r,\z$ are associated invariants given in \eqref{eq-theta-rho}.
Let $H_t=Y_t\wedge\widehat{Y}_t$ be a variation of 
$H=H_0=Y\wedge\widehat{Y}$, so that
\[
\< Y_t,Y_t\>=\< \widehat{Y}_t,\widehat{Y}_t\>=0,
\< Y_t,\widehat{Y}_t\>=-1,~~\Longrightarrow~~\< H_t,H_t\>=-1.
\]
We abbreviate $\frac{\dd}{\dd t}\big\vert_{t=0}$ by a dot.
The only restrictions on the variational vector field $\dot{H}$, 
or equivalently on $\dot{Y},\dot{\widehat{Y}}$, are
($\< \dot{H},H \>=0$)
\begin{equation}
\label{eq-dotY}
\< \dot{Y},Y\>=\< \dot{\widehat{Y}},\widehat{Y} \>
=\< \dot{Y},\widehat{Y} \>+\< Y,\dot{\widehat{Y}} \>=0.
\end{equation}
The first variation of the energy of $H_t$ is:
\[
\frac{\dd}{\dd t}\big\vert_{t=0}E(H_t)
=-2\int_M \< \dot{H},\dd*\dd H \>
=4{\rm i}\int_M \< \dot{H},H_{z\zb} \> \cdot\dd z\wedge\dd\zb.
\]
So $H$ is conformal harmonic iff 
$~\th=0~$ and $~\< \dot{H},H_{z\zb} \>=0,~\forall~\dot{H}.$

First we show {\bf (iii)$\mathbf\Rightarrow$(ii)}. 
Take special variational vector fields
\[
\dot{Y}=0,\quad \dot{\widehat{Y}}=\<\x,\x\> Y+\x.
\]
It is easy to verify that they satisfy \eqref{eq-dotY} by checking
$\< \dot{\widehat{Y}},Y \>=0=\< \dot{\widehat{Y}},\widehat{Y} \>$.
Computation shows
\[
\< \dot{H},H_{z\zb} \> 
=\12\< \dot{\widehat{Y}},\bar{\u}Y_z+{\u}Y_{\zb}+N \> 
=\12\< \dot{\widehat{Y}},\widehat{Y}-\12(\abs{\u}^2+\<\x,\x\>)Y-\x \> 
=-\12\<\x,\x\>.
\]
Since the restriction of the Minkowski metric on the 
M\"obius normal bundle $V^{\bot}$ is positive definite, 
$\dot{E}=0$ implies $\x=0,$
\footnote{
Intuitively, this is because the expression of $\r$ \eqref{rho1} 
contains the term $\<\x,\x\>$, hence $\x$ must vanish if 
the integral of $\r+\bar\r$ is critical.
}
i.e. $\hat{f}$ is on the mean curvature sphere of $f$. 
But there is no bias for $f$ or $\hat{f}$ in the assumptions, 
so these two surfaces should be dual to each other. Hence $f$ 
is also on the mean curvature sphere of $\hat{f}$. 

Next we prove {\bf (ii)$\mathbf\Rightarrow$(i)}. With $\x=0,\th=0$
we have the simplified formulae below:
\begin{align*} 
\widehat{Y}
&=\12\biggl(\abs{\u}^2+\<\x,\x\>\biggr) Y 
+ \bar{\u}Y_z + {\u}Y_{\zb} + N, \\ 
\widehat{Y}_z 
&=\frac{\u}{2}\widehat{Y}  
+ \r\left(Y_z + \frac{\u}{2}Y\right) + 2\e, 
\end{align*}  
where 
\[
\r:=\bar\u_z-2\< \ka,\bar\ka \>, \quad
\e:=D_{\zb}\ka+\frac{\bar\u}{2}\ka.
\]
As in last subsection, $\th:=\u_z-\12\u^2-s=0$ further implies
\[
\r_{\zb}=\bar\u\r+4\<\e,\bar\ka\>
\]
by Gauss equation~\eqref{gauss}, and 
\[
D_{\zb}\e-\frac{\bar\u}{2}\e=D_{\zb}D_{\zb}\ka+\frac{\bar{s}}{2}\ka,
\]
which is real-valued by Codazzi equation~\eqref{codazzi}. 
Also note 
\[
\e_{\zb}=D_{\zb}\e +2\<\e,\bar\e\> Y
-2\<\e,\bar\ka\>\left(Y_z+\frac{\u}{2}Y\right)
\]
by \eqref{eq-structure}. Now the differentiation 
of $\widehat{Y}_z$ can be computed out with the outcome
\begin{align}
\widehat{Y}_{z\zb} 
&= \frac{\u}{2}\widehat{Y}_{\zb}+\frac{\bar\u}{2}\widehat{Y}_z
+ \left( \frac{\u_{\zb}}{2}+\frac{\r}{2}-\frac{\abs{\u}^2}{4} \right)
\widehat{Y} \notag\\
&\qquad +\left( \12\abs{\r}^2+4\<\e,\bar\e\> \right)Y 
+2\left( D_{\zb}\e-\frac{\bar\u}{2}\e \right). \label{yhat-zzbar}
\end{align}
Since $\hat{f}$ is also on the mean curvature sphere of $f$, 
$Y$ is a linear combination of
\[
\{ \widehat{Y},\widehat{Y}_z,\widehat{Y}_{\zb},\widehat{Y}_{z\zb} \}.
\]
The $\widehat{Y}_{z\zb}$-component of $Y$ is not zero. (Otherwise
$Y$ is a combination of
$\widehat{Y},\widehat{Y}_z,\widehat{Y}_{\zb}$, 
hence $\< Y,\widehat{Y} \>=0$, a contradiction.)
So $\widehat{Y}_{z\zb}$, as well as $D_{\zb}\e-\frac{\bar\u}{2}\e$, 
is contained in
$~\text{Span}
\{ \widehat{Y},\widehat{Y}_z,\widehat{Y}_{\zb},Y \}$.
By the expressions of $\widehat{Y},\widehat{Y}_z$, this is true only if
$
0=D_{\zb}\e-\frac{\bar\u}{2}\e=D_{\zb}D_{\zb}\ka+\frac{\bar{s}}{2}\ka,
$
i.e. $f$ is Willmore. Again by the duality between $f$ and $\hat{f}$
we know $\hat{f}$ is also Willmore. The assumptions directly imply
that they form adjoint transform to each other.

Finally one should verify {\bf (i)$\mathbf\Rightarrow$(iii)}. 
This case $\th=0,~\x=0,~D_{\zb}\e-\frac{\bar\u}{2}\e=0$,
and \eqref{yhat-zzbar} takes the following form:
\[
\widehat{Y}_{z\zb} 
= \frac{\u}{2}\widehat{Y}_{\zb}+\frac{\bar\u}{2}\widehat{Y}_z
+ \left( \frac{\u_{\zb}}{2}+\frac{\r}{2}-\frac{\abs{\u}^2}{4} \right)
\widehat{Y} +\left( \12\abs{\r}^2+4\<\e,\bar\e\> \right)Y. 
\]
We compute $\< \dot{H},H_{z\zb} \>$ for arbitrary $\dot{H}$, 
or equivalently, for any variational vector fields 
$\dot{Y},\dot{\widehat{Y}}$. 
Invoking the restrictions \eqref{eq-dotY}, there follows
\begin{align*}
\< \dot{H},H_{z\zb} \>
&= \< \dot{Y},\widehat{Y}_{z\zb} 
-\frac{\u}{2}\widehat{Y}_{\zb}-\frac{\bar\u}{2}\widehat{Y}_z
+\left( \frac{\abs{\u}^2}{4}-\<\ka,\bar\ka\> \right)\widehat{Y} \> \\
&\qquad 
+\< \dot{\widehat{Y}}, \frac{\u}{2}Y_{\zb}+\frac{\bar\u}{2}Y_z
+\12 N-\left( \< Y,\widehat{Y}_{z\zb} \>+\<\ka,\bar\ka\> \right)Y \> \\
&= \< \dot{Y},
\left( \frac{\u_{\zb}}{2}+\frac{\r}{2}-\<\ka,\bar\ka\> \right)\widehat{Y}
+\left( \12\abs{\r}^2+4\<\e,\bar\e\> \right)Y \> \\
&\qquad 
+\< \dot{\widehat{Y}},\12\widehat{Y}-\frac{\abs{\u}^2}{4}Y
+\left( \frac{\u_{\zb}}{2}+\frac{\r}{2}+\frac{\abs{\u}^2}{4}
-\<\ka,\bar\ka\> \right)Y \> \\
&= 0.
\end{align*}
So $H$ is harmonic. This completes our proof.
\end{proof}
\begin{remark}
Indeed, our discussion above provides another 
proof to the Duality Theorem~\ref{thm-duality}.
\end{remark}
\begin{remark}
As Burstall pointed out to the author, one striking feature of the 
theorem above is that the condition (iii) is a collection of first 
and second order conditions, yet they force a fourth order equation 
(condition (i), that $f,\hat{f}$ being Willmore surfaces).
\end{remark}


\renewcommand{\thetheorem}{A-\arabic{theorem}}
\setcounter{theorem}{0}  
\section*{Appendix}

In this part we want to discuss the relationship of left/right touching
between 2-planes in the setup of quaternions, then point out the
connection with the general notion of touch and co-touch. Our starting
point is the following lemma. 
\begin{lemma} 
[the Fundamental Lemma in \cite{quater}; also Lemma~6 in \cite{Bohle}] 
For every oriented real subspace $U \subset \H$ of dimension 2 there are 
unique vectors $N,R$ satisfying $N^2=R^2=-1$ and 
\begin{equation} 
\label{eq-normal} 
U = \{x \in \H ~|~Nx=-xR~\}.
\end{equation} 
Conversely, every pair of vectors $N$ and $R$ satisfying $N^2=R^2=-1$ 
defines, via \eqref{eq-normal}, an oriented 2-plane. 
\end{lemma} 

$N$ and $R$ are called the 
\emph{left} and \emph{right normal vector} of $U$ respectively, 
though in general they are \emph{not} orthogonal to $U$. 
For an oriented immersed surface $M$ in $\H$, we can define 
the pair $\{N,R\}$ similarly, which might be identified with 
the usual Gauss map of $M$ in $\R^4$. 

\begin{definition} 
Let $U_i$ be oriented 2-plane with $N_i$ and $R_i$ as their left and right 
normal vectors respectively, $i=1,2.$ Then 
\begin{enumerate} 
\item $U_1$ and $U_2$ \emph{touch} each other \emph{from left (right)} if 
$N_1=N_2~(R_1=R_2)$.
\item $U_1$ and $U_2$ \emph{co-touch} each other \emph{from left (right)} if 
$N_1=-N_2~(R_1=-R_2)$.
\end{enumerate} 
Similarly we can define \emph{(co-)touch} of two conformal immersions 
at their intersection point.
When touch is both from left and right, these two 
immersions are tangent at the intersection point with the same 
induced orientation. So left/right touch may be viewed as the 
generalization of tangency. 
\end{definition}

Given two oriented 2-planes in an oriented 4-space, 
This definition seems algebraic and depending on the way 
in which we identify $\R^4\simeq\H$. Yet by the following two lemmas, 
we find they are well-defined geometric notions (depending only on 
different choices of orientations). 

\begin{lemma} 
\label{lem-isometry}
Every orientation preserving linear isometry of $\H$ is of the form 
$x\in\H ~\mapsto~ \u x\l. $ Here $\u,\l\in\S^3$ are unit quaternions.
\end{lemma}

\begin{lemma}
\mbox{}
\begin{itemize} 
   \item[(i)] 
Every orientation preserving linear isometry of $\H$ leaves the 
relationship of \emph{left (co-)touch} and \emph{right (co-)touch} 
invariant. 
   \item[(ii)] 
Every orientation reversing linear isometry of $\H$ preserves 
the property of \emph{touch} and \emph{co-touch}, but interchanges 
between \emph{left} and \emph{right}. 
   \item[(iii)] 
Suppose $U_1$ touches (co-touches) $U_2$ from left/right. 
Then $U_1$ with opposite orientation co-touches (touches) $U_2$ 
from left/right respectively. 
\end{itemize}
\end{lemma} 

The second lemma is easy to obtain as a corollary of the first one, 
which is a well-known fact and we omit both proofs at here. 
They tell us that the difference between \emph{left} and \emph{right} 
is due to the orientation induced by the identification $\R^4=\H$, 
hence not essential. The proposition below confirms this observation, 
and unifies two different definitions of touch and co-touch. 

\begin{proposition} 
Let $U$ and $\widehat{U}$ be a pair of oriented 2-dim subspaces in 
$\R^4.$ They (co-)touch each other as contact elements if, and only if, 
they (co-)touch each other from left or right. (Whether it is from left 
or right depends on the orientation induced by the identification 
$\R^4\simeq\H$.) 
\end{proposition} 
\begin{proof} 
Equipped $U,\widehat{U}$ with oriented orthonormal basis 
$\{\a,\b\}$ and $\{\hat\a,\hat\b\}$ respectively. 
Regarding $U$ as a conformally embedded submanifold of $\R^4\subset\S^4$, 
we fix a lift $U\subset \R^4 \to \R^{5,1}$ as 
\[ 
v\in U ~\mapsto~ \left( \12(1+\abs{v}^2),\12(1-\abs{v}^2),v \right), 
\] 
which projects down to $\P(\mathcal{L})$. The image of $0\in\R^4$ is 
$Y=(\12,\12,0,\dotsc,0)$. The induced contact element at $0$ is given by 
$\Sg=\{Y,Y_1,Y_2\}$, where
\[ 
Y_1=(0,0,\a),\quad Y_2=(0,0,\b). 
\] 
We have similar representation $\widehat\Sg=\{Y,\widehat Y_1,\widehat Y_2\}$ 
for $\widehat{U}$. Now consider the invariant $\underline\r$ associated with 
$\Sg,\widehat\Sg$. By definition \eqref{theta-rho}, 
\[ 
\underline\r=\12\<\a-{\rm i}\b,\hat\a-{\rm i}\hat\b \>=0 
~~\Longleftrightarrow~~ 
\left\{ 
\begin{array}{ccr} 
\<\a,\hat\a\> &=& \<\b,\hat\b\>,\\ 
\<\a,\hat\b\> &=& -\<\b,\hat\a\>. 
\end{array} 
\right. 
\] 
We define a complex structure $J$ on $\R^4$ via 
$J\{\a,\b,\hat\a,\hat\b\}=\{\b,-\a,\hat\b,-\hat\a\}. $
In this term the condition of touch holds if, and only if, 
there is $J$ satisfying the formula above and compatible 
with the Euclidean metric. By Lemma~\ref{lem-isometry}, 
it is easy to show that such a complex structure 
must be of the form $\a\mapsto N\a$ or $\a\mapsto\a R$, 
where $N,R\in\H, N^2=R^2=-1.$ This implies our conclusion 
on touch. For co-touch the similar argument applies. 
\end{proof}

\vspace{.4cm}
\noindent {\bf Acknowledgment} \\
This work is part of the Doctoral dissertation of 
the author at TU-Berlin. The author would like to thank his advisor
Prof. U. Pinkall for the inspiring discussion, and thank TU-Berlin for 
providing a special scholarship for a period of four years.


\vspace{0.4cm}

Xiang Ma, Fachbereich Mathematik, 
Technische Universit\"{a}t Berlin, 
Str. des 17. Juni 136, D-10623 Berlin, Germany\\

\textit{Email address}: \textsf{ma@math.tu-berlin.de}

\end{document}